\documentclass[11pt]{amsart}
\usepackage{amssymb, latexsym, amsmath, amsfonts, tikz, tikz-cd}
\usepackage{hyperref, enumitem, fancyhdr}

\newtheorem{thm}{Theorem}[section]

\theoremstyle{definition}

\theoremstyle{remark}

\numberwithin{equation}{section}
\theoremstyle{remark}

\newcommand{\mbb}{\mathbb}
\newcommand{\ra}{\rightarrow}

\newcommand{\pa}{\partial}
\newcommand{\ov}{\overline}
\newcommand{\sm}{\setminus}
\newcommand{\ep}{\epsilon}
\newcommand{\no}{\noindent}
\newcommand{\al}{\alpha}

\newcommand{\cal}{\mathcal}
\newcommand{\ti}{\tilde}
\newcommand{\la}{\lambda}

\newcommand{\ga}{\gamma}

\newcommand{\be}{\beta}

\begin{document}
\title{Two remarks on the Poincar\'{e} metric on a singular Riemann surface foliation}
\author{Sahil Gehlawat and Kaushal Verma}

\address{SG: Department of Mathematics, Indian Institute of Science, Bangalore 560 012, India}
\email{sahilg@iisc.ac.in}

\address{KV: Department of Mathematics, Indian Institute of Science, Bangalore 560 012, India}
\email{kverma@iisc.ac.in}

\begin{abstract} 
 Let $\mathcal F$ be a smooth Riemann surface foliation on $M \sm E$, where $M$ is a complex manifold and $E \subset M$ is a closed set. Fix a hermitian metric $g$ on $M \sm E$ and assume that all leaves of $\mathcal F$ are hyperbolic. For each leaf $L \subset \mathcal F$, the ratio of $g | L$, the restriction of $g$ to $L$, and the Poincar\'{e} metric $\la_L$ on $L$ defines a positive function $\eta$ that is known to be continuous on $M \sm E$ under suitable conditions on $M, E$. For a domain $U \subset M$, we consider $\mathcal F_U$, the restriction of $\mathcal F$ to $U$ and the corresponding positive function $\eta_U$ by considering the ratio of $g$ and the Poincar\'{e} metric on the leaves of $\mathcal F_U$. First, we study the variation of $\eta_U$ as $U$ varies in the Hausdorff sense motivated by the work of Lins Neto--Martins. Secondly, Minda had shown the existence of a domain Bloch constant for a hyperbolic Riemann surface $S$, which in other words shows that every holomorphic map from the unit disc into $S$, whose distortion at the origin is bounded below, must be locally injective in some hyperbolic ball of uniform radius. We show how to deduce a version of this Bloch constant for $\mathcal F$.
\end{abstract}  

\maketitle 

\section{Introduction}

\noindent To provide a context for this note, let $M$ be a complex manifold of dimension $n$, $E \subset M$ a closed subset and $\mathcal F$ a smooth lamination of $M \sm E$ by Riemann surfaces. Recall that in this case, $\mathcal F$ is defined by an atlas of flow-boxes $(U_{\al}, \phi_{\al})$ where $\phi_{\alpha} : U_{\alpha} \ra \mbb D \times \mbb D^{n-1}$ ($\mbb D \subset \mbb C$ is the open unit disc) and the transition maps 
$\phi_{\al\be}$ are smooth in $(x, y) \in \mbb D \times \mbb D^{n-1}$ and are of the form
\[
\phi_{\al \be} = \phi_{\al \be}(x, y) = \phi_{\al} \circ \phi^{-1}_{\be}(x, y) =  (A(x, y), B(y))
\]
with $A$ holomorphic in $x$. The leaves of $\mathcal F$, which are locally given by $\phi^{-1}_{\al}(\{y = \mbox{constant} \})$, are Riemann surfaces and the one passing through a given $p \in M \sm E$ will be denoted by $L_p$. We will assume that each leaf is hyperbolic and refer to $(M, E, \mathcal F)$ as a hyperbolic Riemann surface lamination.
 Fix a hermitian metric $g$ on $M \sm E$ and denote the length of a tangent vector $v$  by $\vert v \vert_{g}$. Consider the space $\mathcal O(\mbb D, \mathcal F)$ of all holomorphic maps $f$ from the unit disc $\mbb D \subset \mbb C$ with values in a leaf of $\mathcal F$, i.e., $f(\mbb D) \subset L$ for some leaf $L \subset \mathcal F$. Pick $f \in \mathcal O(\mbb D, \mathcal F)$ and let $f(0) = z$ so that $f(\mathbb D) \subset L_z$. Let $\pi_z : \mbb D  \ra L_z$  be the universal covering map chosen such that $\pi_z(0) = z$. Then $f$ admits a unique lift $\ti f : \mbb D \ra \mbb D$ that fixes the origin. The Schwarz lemma implies that $\vert \ti f'(0) \vert \le 1$ and hence $\vert f'(0) \vert_g \le \vert \pi'_z(0) \vert_g$. Therefore, $\eta : M \sm E \ra (0, \infty)$ given by
\[
\eta(z) = \sup \left \{ \vert f'(0) \vert_g : f \in  \mathcal O(\mbb D, \mathcal F), f(0) = z \right\}
\]
is well-defined. Note that $\eta > 0$ everywhere by considering the holomorphic map that sends $\mbb D$ into a local leaf $\phi^{-1}_{\al}(\{y = \mbox{constant} \})$, and that its extremal value is achieved by $\pi_z$, which evidently belongs to $\mathcal O(\mbb D, \mathcal F)$ . Thus, $\eta(z) = \vert \pi'_z(0) \vert_g$. Let $\la_{ L}$ denote the Poincar\'{e} metric on a leaf $L \subset \mathcal L$ and for a vector $v$ that is tangent to $L$ at $z$, let $\vert v \vert_{L}$ denote its Poincar\'{e} length. The extremal property of the Kobayashi (which equals the Poincar\'{e}) metric shows that $\vert v \vert_g = \eta(z) \vert v \vert_{L}$. In other words, the restriction of $g$ to a leaf $L$ differs from $\la_ L$ by a multiplicative constant which is captured by the positive function $\eta$. This property will be written symbolically as $g = \eta \cdot \la_{L}$ and can be understood to mean that the restriction of $g/\eta$ to a leaf $L$ induces the Poincar\'{e} metric $\la_{L}$ on it.

\medskip

Now, for a fixed leaf $L$, $\la_{L}$ is a positive real-analytic function on it while $g$ is smooth, and hence $\eta$ is also smooth along leaves. However, the regularity of $\eta$ along transverse directions is not a priori clear. In fact, it is related to the variation of $\la_{L}$ along such directions. This question has been studied under various hypotheses on $M, E$ and 
$\mathcal F$, and we refer the reader to Verjovsky \cite{V}, Lins Neto \cite{N1, N2}, Candel \cite{Ca} in which the continuity of $\la_ {L}$ on $M \sm E$ is established. In addition, recent work of Dinh--Nguy\^{e}n--Sibony \cite{DNS1} shows that $\la_{L}$ is in fact H\"{o}lder continuous when $M$ is compact and $E = \emptyset$; a stronger estimate on the modulus of continuity 
of $\la_{L}$ can be found in \cite{DNS2} under suitable hypotheses on $M, E$ and $\mathcal F$. Related work on this theme can be found in \cite{CLS}, \cite{CG}, and \cite{FS}.

\medskip

A different continuity property of $\eta$ was proved in \cite{NM}. To describe it, let $\mathcal F_U$ denote the restriction of $\mathcal F$ to a domain $U \subset M$. For $p \in M \sm E$, let $L_{p, U}$ denote the connected component of $L_p \cap U$ that contains $p$. If it is non-empty, $L_{p, U}$ is then the leaf of $\mathcal F_U$ that contains $p$. Let $\eta_{U} : U \ra (0, \infty)$ be the positive function corresponding to the leaves of $\mathcal F_U$. Note that $\eta_{U}$ is defined using holomorphic maps from $\mbb D$ into $L_{p, U}$ as $p$ varies in $U$. The following facts were established in \cite{NM}. First, that $\eta_{U}$ is monotone as a function of $U$, i.e., for domains $U \subset V$, $\eta_{U}(p) \le \eta_{V}(p)$ for all $p \in U$. Second, if $\{U_n\}$ is an increasing exhaustion of $U$, then $\eta_{U_n} \ra \eta_{U}$ pointwise on $U \sm E$, and uniformly on compact subsets of $U \sm E$ if the functions $\eta_{U_n}$ and 
$\eta_{U}$ are known to be continuous. 

\medskip

To strengthen this, let $d$ be the distance on $M$ induced by $g$, and for 
$S \subset M$ and $\ep > 0$, let $S_{\ep}$ be the $\ep$-thickening of $S$ with distances being measured using $d$. Recall that the Hausdorff distance $\mathcal H(A, B)$ between compact sets $A, B \subset M$ is the infimum of all $\ep >0$ such that $A \subset B_{\ep}$ and $B \subset A_{\ep}$. For bounded domains $U, V \subset M$, the prescription $\rho(U, V) = \mathcal H(\ov U, \ov V) + \mathcal H(\pa U, \pa V)$ defines a metric (see \cite{Bo}) on the collection of all  bounded open subsets of $X$ with the property that if $\rho(U, U_n) \ra 0$, then every compact subset of $U$ is eventually contained in $U_n$, and every neighbourhood of $\ov U$ contains all the $U_n$'s eventually.

\begin{thm}
Let $(M, E, \mathcal F)$ be a hyperbolic Riemann surface lamination as above. 
\begin{itemize}
\item[(i)] Suppose that $U \subset M$ is a bounded taut domain and $\{U_n\}$ is a sequence of bounded domains in $M$ such that $\rho(U, U_n) \ra 0$. Then, $\eta_{U_n}$ converges to $\eta_{U}$ pointwise on $U \sm E$.
\item[(ii)] In addition, assume that $E$ is discrete, $\eta_{U}$ is continuous, and there is a taut domain $V$ containing $\ov U$. Then, $\eta_{U_n}$ converges to $\eta_{U}$ uniformly on compact subsets of $U \sm E$.
\end{itemize}
\end{thm} 

Since the functions $\eta_U, \eta_{U_n}$ are a measure of the uniformizations of the leaves, the crux of the proof of these statements is to study the compactness of a family of uniformizations. In (i), which is a pointwise statement, the leaf $L_p$ does not change but only the open connected subsets $L_{p, U}$ are allowed to vary. On the other hand, the situation considered in (ii) allows even the leaves to vary. 

\medskip

The other notion that we consider is motivated by Minda's definition of a domain Bloch constant -- see \cite{Mi1}. To recall this briefly, let $S$ be a hyperbolic Riemann surface and $\la_S(z) \vert dz \vert$ the Poincar\'{e} metric on it. For $1 \le m \le \infty$, consider the family $\mathcal O_m(\mbb D, S)$  As examples, $\mathcal O_1(\mbb D, S) $ is precisely the family of all non-constant holomorphic maps, while $\mathcal O_{\infty}(\mbb D, S)$ is to be understood as the family of locally injective holomorphic maps $f : \mbb D \ra S$, for roots with infinite multiplicity can only exist for constant maps. 

\medskip

Let $d_h(z, w)$ denote the distance between $z, w \in \mbb D$ in the Poincar\'{e} metric and  let $B_h(z, \delta) \subset \mbb D$ be the open ball of radius $\delta > 0$ centered at $z$ with respect to $d_h$. Now given $f \in \mathcal O_m(\mbb D, S)$ and $z \in \mbb D$, let
\[
t(z, f) = \sup \{ r \in (0, +\infty): f \;\mbox{is locally injective in} \; B_h(z, r) \}
\] 
where we set $t(z, f) = 0$ if $f$ is not locally injective in any neighbourhood of $z$. The parameter
\[
t(f) = \sup\{ t(z, f): z \in \mbb D\}
\]  
is then a measure of the largest hyperbolic ball in $\mbb D$ on which $f$ is locally injective. The {\it distortion function} defined by
\[
f^S_{\mbb D}(p) = \left. \frac{f^{\ast}(\la_S(z) \vert dz \vert)}{\la_{\mbb D}(z) \vert dz \vert}  \right|_{z = p}
\]
measures the local distortion induced by $f$ at $p$ measured with respect to the Poincar\'{e} metrics on $S$ and $\mbb D$. It is known (see \cite{Mi2}) that the distortion is at most $1$ and equality holds precisely when $f : \mbb D \ra S$ is a holomorphic covering. For $\al \in (0, 1]$,  
\begin{equation}
\mathcal T_{m, S}(\al) = \inf \{ t(f): f \in \mathcal O_m(\mbb D, S) \;\mbox{and} \;f^S(0) \ge \al \}
\end{equation}
is called the {\it domain Bloch constant}, which if positive, would show that every $f \in \mathcal O_m(\mbb D, S)$ with distortion at least $\al$ at the origin is locally injective in {\it some} hyperbolic ball in $\mbb D$ of uniform radius. Note that $\mathcal T_{m, S}$ is an increasing function of both $m$ and $\al$ and $\mathcal T_{\infty, S}(\al) = +\infty$ for all $\al \in (0, 1]$. Minda \cite{Mi1} showed that for $1 \le m < \infty$,
\begin{equation}
\mathcal T_{m, S}(\al) \ge \log \left( \frac{1 + \tau_m (\al^{1/(m+1)})}{1 - \tau_m (\al^{1/(m+1)})} \right)
\end{equation}
where  
\[
\tau_m(\al) = \al \left( \frac{m(m+2)}{2(m+1)^2 - (m^2 + 2m + 2)\al^2 + 2(m+1)\left( \al^4 - (m^2 + 2m + 2)\al^2 + (m+1)^2  \right)^{1/2}} \right)^{1/2}
\]
and as can be checked, is an increasing function from the closed interval $[0,1]$ to itself that fixes both $0$ and $1$.

\medskip

Since these considerations show the existence of a domain Bloch constant for holomorphic maps into a fixed hyperbolic Riemann surface, it seems natural to ask whether there is a foliated version of this notion. It turns out that there is. To describe it, let $(M, E, \mathcal F)$ be a Riemann surface lamination and recall the notion of an {\it ultrahyperbolic metric on $\mathcal F$ of curvature bounded by $-a < 0$} for some $a > 0$. Following \cite{NM}, \cite{N2}, a continuous hermitian metric $g$ on $M \sm E$ is such a metric if it satisfies the following conditions: 

\medskip

\no (i) For a leaf $L \subset \mathcal F$ and a local holomorphic coordinate $(z, U)$ on it, $g | L$, the restriction of $g$ to $L$ is of the form $g | L  = f(z) \vert dz \vert^2$ in $U$, where $f$ is continuous and non-negative on $U$  and $f^{-1}(0)$ is discrete in $U$, and

\medskip

\no (ii) For every $z \in U$ for which $f(z) \not=0$, there is a positive $C^2$-function $h$ near $z$ such that $h(z) = f(z)$, $h \le f$ near $z$, and the Gaussian curvature of $h(z) \vert dz \vert^2$ is at most $-a < 0$. 

\medskip

\no Let $g$ be such a metric on $M \sm E$. As before, for $1 \le m \le \infty$, let $\mathcal O_m(\mbb D, \mathcal F)$ denote the family of all non-constant holomorphic maps from $\mbb D$ with values in a leaf of $\mathcal F$ with the property that for any $q \in f(\mbb D)$, each root of $f(z) = q$ is either simple or else has multiplicity at least $m+1$. For a given $f \in \mathcal O_m(\mbb D, \mathcal F)$, suppose that $f$ takes values in $L \subset \cal F$. Then 
\[
f^{\mathcal F, g}_{\mbb D}(p) = \left. \frac{f^{\ast}(g | L)}{\la_{\mbb D}(z) \vert dz \vert}  \right|_{z = p}
\]
defines a distortion function of such an $f$ with respect to $g$. As $f$ varies in $\mathcal O_m(\mbb D, \mathcal F)$, we get a family of such distortion functions. With $t(z, f)$ and $t(f)$ defined as before, and $\al > 0$, set 
\[
\mathcal T^g_{m, \cal F}(\al) = \inf \{ t(f) : f \in  \mathcal O_m(\mbb D, \mathcal F), f^{\mathcal F, g}_{\mbb D}(0) \ge \al \}.
\]

\begin{thm}
Let $(M, E, \cal F)$ be a Riemann surface foliation and suppose that $g$ is an ultrahyperbolic metric on $\cal F$ with curvature bounded above by $-a^2 < 0$ for some $a > 0$. Then
\[
\mathcal T^g_{m, \cal F}(\al) \ge  \log \left( \frac{1 + \tau_m ((a \al)^{1/(m+1)})}{1 - \tau_m ((a \al)^{1/(m+1)})} \right)
\]
for $\al \in (0, 1/a]$. 
\end{thm}
 
Thus, every $f \in \mathcal O_m(\mbb D, \mathcal F)$ whose distortion at the origin is at least $\al \in (0, 1/a]$ is necessarily locally injective in some hyperbolic ball of uniform radius --  a lower bound on this radius is given above. The final section contains miscellaneous remarks and examples pertaining to the theorems presented above.


\section{Proof of Theorem 1.1}

\no For {\it (i)}, let $ p \in U \sm E$. Then $p \in U_n$ eventually. Let $\al : \mbb D \ra L_{p, U}$ be a uniformizing map with $\al(0) = p$ so that $\eta_{U}(p) = \vert \al'(0) \vert_g$. For 
$\ep > 0$, let $\psi : \mbb D \ra \mbb D$ be defined by $\psi(z) = (1-\ep)z$, and note that $\al \circ \psi(\mbb D)$ is compactly contained in $U$. Hence, $\al \circ \psi(\mbb D) \subset U_n$ for all large $n$ since  $\rho(U, U_n) \ra 0$. Further, as $\al \circ \psi (\mbb D)$ is connected and contained in $L_p$, it follows that $\al \circ \psi(\mbb D)$ is in fact contained in $L_{p, U_n}$. The map $\al \circ \psi$ is then a candidate in the family that defines $\eta_{U_n}$ and therefore
\[
\eta_{U_n}(p) \ge \vert (\al \circ \psi)'(0) \vert_g = (1-\ep) \vert \al'(0) \vert_g = (1-\ep) \eta_{U}(p)
\]
for all large $n$. This shows that $\liminf \eta_{U_n}(p) \ge \eta_{U}(p)$.

\medskip

Now let $\{\eta_{U_{n_j}} \}$ be a convergent subsequence, and let $\al_{n_j} : \mbb D \ra L_{p, U_{n_j}}$ be a family of uniformization maps such that $\al_{n_j}(0) = p$. Then, 
\[
\eta_{U_{n_j}}(p) = \vert \al'_{n_j}(0) \vert_g 
\]
and since each $\al_{n_j}$ takes values in the fixed hyperbolic leaf $L_p$, the family of uniformizations is normal. A further subsequence, which will still be indexed by $n_j$ for the sake of brevity, then converges to a holomorphic limit $\ti \al : \mbb D \ra L_p$ with $\ti \al (0) = p$.

\medskip

The claim now is that $\ti \al(\mbb D) \subset L_{p, U}$. To see this, suppose to the contrary that $q = \ti \al(z_0) \notin \ov U$. Since $\rho(U, U_n) \ra 0$, $q \notin \ov U_n$ for all large $n$. On the other hand, the fact that $\al_{n_j}(z_0) \ra \ti \al(z_0) = q$ implies that $\al_{n_j}(z_0) \notin U_{n_j}$ for all large $n_j$ -- this is a contradiction since $\al_{n_j}(\mbb D) \subset U_{n_j}$. Hence $\ti \al(\mbb D) \subset \ov U$. But as $U$ is taut and $\ti \al(0) = p \in U$, it follows that $\ti \al(\mbb D) \subset U$ and hence a subset of $L_p$. Thus, $\ti \al$ is a candidate in the family that defines $\eta_{g, U}$ and hence
\[
\eta_{U_{n_j}}(p) = \vert \al'_{n_j}(0) \vert_g \ra \vert \ti \al'(0) \vert_g  \le \eta_{U}(p). 
\]
This shows that $\limsup \eta_{U_n}(p) \le \eta_{U}(p)$ -- this completes the proof of {\it (i)}.

\medskip

For {\it (ii)}, suppose that there is a compact set $K \subset U \sm E$ and a sequence $\{p_n\} \in K$ such that
\[
\vert \eta_{U_n}(p_n) - \eta_U(p_n) \vert > \ep
\]
for some $\ep > 0$. Assume that $p_n \ra p$ after passing to a subsequence. By the continuity of $\eta_U$, $\vert \eta_U(p_n) - \eta_U(p) \vert < \ep/2$
and hence
\begin{equation}
\vert \eta_{U_n}(p_n) - \eta_U(p) \vert > \ep/2
\end{equation}
for large $n$. Let $\al_n : \mbb D \ra L_{p_n, U_n} \subset U_n$ be a uniformization map with $\al_{p_n}(0) = p_n$. Since $\rho(U, U_n) \ra 0$, $V$ eventually contains all the $U_n$'s and hence the family $\{\al_n\}$ is normal. By passing to a subsequence, let $\ti \al : \mbb D \ra \ov V$ be a holomorphic limit with $\ti \al(0) = p$. To complete the proof, it suffices to show that $\ti \al : \mbb D \ra L_{p, U}$ is a uniformization. Indeed, if this were the case, then $\eta_U(p) = \vert \ti \al'(0) \vert_g$ which would then imply that
\[
\eta_{U_n}(p_n) = \vert \al'_n(0) \vert_g \ra \vert \ti \al'(0) \vert_g = \eta_U(p)
\]
and this contradicts (2.1). That $\ti \al$ uniformizes $L_{p, U}$  is a consequence of the following observations that are inspired by the proof of Proposition $3$ in \cite{NM}:

\medskip

\no (a) Since $V$ is taut and $\ti \al(0) = p \in V$, it follows that $\ti \al (\mbb D) \subset V$. Now suppose that $\ti \al(a) \in V \sm \ov U$ for some $a \in \mbb D$. Then $\al_n(a) \in V \sm \ov U$ for all large $n$ and since $\rho(U, U_n) \ra 0$, it follows that $\al_n(0) \notin \ov U_n$ for all large $n$. This is a contradiction and hence $\ti \al : \mbb D \ra \ov U$. Again, the tautness of $U$ implies that $\ti \al : \mbb D \ra U$.

\medskip

\no (b) To see that $\ti \al(\mbb D) \subset L_{p, U} \cup E$, note that there are two possibilities at this stage. If $\ti \al$ is constant, then $\ti \al(\mbb D) = \ti \al(0) = p \in L_{p, U}$. If $\ti \al$ is non-constant, note that 
\[
N = \{z \in \mbb D: \ti \al(z) \notin E \}
\]
is non-empty, for $\ti \al(0) = p \in U \sm E$ shows that the origin is contained in $N$. For $w_0 \in N$, let
\[
N_{w_0} = \{ w \in N : L_{\ti \al(w)} = L_{\ti \al(w_0)}\}.
\]
To see that $N_{w_0}$ is open in $\mbb D$, pick $w_1 \in N_{w_0}$ and a chart $ \phi = (x, y) : W \ra \mbb D \times \mbb D^{n-1}$ around $\ti \al(w_1) = q$ such that $(x, y)(q) = 0 \in \mbb C^n$ and the leaves of $\mathcal F_W$ are described by $y = \mbox{constant}$. Let $V$ be the connected component of $\ti \al^{-1}(W) \subset \mbb D$ containing $w_1$. Now, the claim is that $y(\ti \al(V)) = 0$. If this is not the case, then $y(\ti \al (V)) \not= \{0\}$. Since $y(\ti \al(w_1)) = 0$, it follows that $y \circ \ti \al : V \ra \mbb D^{n-1}$ is a non-constant map which implies that its derivative $d(y \circ \ti \al)(w_2) \not= 0$ for some $w_2 \in V$ near $w_1$. Since $\al_n$ converges to $\ti \al$, $d(y \circ \al_n)(w_2) \not= 0$ for $n$ large. This means that $\al_n(\mbb D)$ is not contained in a leaf of $\mathcal F$. This contradiction shows that $y(\ti \al(V)) = 0$. Hence every $w \in V$ has the property that the leaf of $\mathcal F_W$ containing $\ti \al(w)$ coincides with $y = 0$ which is also the leaf passing through $w_1$. This shows that $N_{w_0}$ is open.

\medskip

Since $E$ is discrete, $\ti \al^{-1}(E) \subset \mbb D$ is also discrete since $\ti \al$ is assumed to be non-constant. Therefore, $\mbb D \sm \ti \al^{-1}(E)$ is connected. Since $N_{w_0}$ is open for each $w_0 \in N$, it follows that $\mbb D \sm \ti \al^{-1}(E)$ coincides with $N_{\ti z}$ for some $\ti z \in N$. Hence, $L_{\ti \al(z)} = L_{\ti \al(\ti z)}$ for all $z \in \mbb D \sm \ti \al^{-1}(E)$ and this implies that $\ti \al (z) \in L_{\ti \al(z)} = L_{\ti \al(\ti z)}$ for all $z \in  \mbb D \sm \ti \al^{-1}(E)$. By continuity, $\ti \al(\mbb D) \subset L_{\ti \al(\ti z)} \cup E$. But $\ti \al(0) = p \in L_p$ since $p \notin E$. This shows that $L_{\ti \al(\ti z)} = L_p$ and hence that $\ti \al(\mbb D) \subset L_p \cup E$.

\medskip 

\no (c) We now show that $\ti \al : \mbb D \sm \ti \al^{-1}(E) \ra L_p$ is a covering map. This will be done in two steps -- first, we show that $\ti \al$ is an immersion on $\mbb D \sm \ti \al^{-1}(E)$ and second, that it has the path lifting property on $L_p$. 
 
\medskip

For the first step, pick $z_0 \in \mbb D \sm \ti \al^{-1}(E)$ so that $\ti \al(z_0) \notin E$. Let $W$ be a flow box around $\ti \al(z_0)$ small enough so that $W$ is compactly contained in $U$. Since $\rho(U, U_n) \ra 0$, $W \subset U_n$ for large $n$. There is a chart $\phi = (x, y) : W \ra \mbb D \times \mbb D^{n-1}$ with $(x, y)(\ti \al(z_0)) = 0 \in \mbb C^n$ and the leaves of $\mathcal F_W$ are described by $y = \mbox{constant}$. Since $\al_n(z_0) \ra \ti \al(z_0)$, $\al_n(z_0) \in W$ for large $n$. Let $s_n = \al_n(z_0)$ and let $\psi_n$ be a family of holomorphic automorphisms of $\mbb D$ with $\psi_n(0) = \phi(s_n)$; here, we consider the unique local leaf, which is a copy of $\mbb D$, that contains $\phi(s_n)$. Then, $\phi^{-1} \circ \psi_n : \mbb D \ra W$ maps $\mbb D$ into the local leaf containing $s_n$. Now,
\[
\vert (\phi^{-1} \circ \psi_n)'(0) \vert_g  = \vert (\phi^{-1})'(\phi(s_n)) \circ \psi'_n(0) \vert_g = \vert  (\phi^{-1})'(\phi(s_n)) \vert_g (1 - \vert \phi(s_n) \vert^2) \ge c > 0
\]
for some uniform $c > 0$, since $\vert \phi(s_n) \vert \ra 0$ and $\phi^{-1}$ is a smooth diffeomorphism. By the extremal property of $\eta_{U_n}$,
\begin{equation}
\vert \al'_n(z_0) \vert_g \ge \vert (\phi^{-1} \circ \psi_n)'(0) \vert_g \ge c > 0
\end{equation}
and hence $\vert \ti \al(z_0) \vert_g \ge c > 0$. Thus, $\ti \al$ is an immersion near $z_0$, and in particular near the origin since $\ti \al(0) \notin E$.

\medskip

In fact, more can be said at this stage. Observe that Rouche's theorem shows that if $G \subset \mbb C$ is a domain with $0 \in G$ and $f : G \ra \mbb D$ is holomorphic with $f'(0) \not= 0$, then there is a uniform $\delta > 0$ such that the euclidean disc $B(f(0), \delta \vert f'(0) \vert/2) \subset f(B(0, \delta))$. With $\phi$ a chart as above near $\ti \al(0) = p$,  let $\ti G \subset \mbb D$, $0 \in \ti G$ be such that $\phi \circ \ti \al : \ti G \ra \mbb D$. Let $G \subset \ti G$ be a compactly contained neighbourhood of $0$. Then, each $\phi \circ \al_n : G \ra \mbb D$, since $\phi \circ \al_n \ra \phi \circ \ti \al$ uniformly on compact subsets of $G$. Therefore, 
there is a uniform $\delta > 0$ such that
\[
B(\phi \circ \al_n(0), \delta  \vert (\phi \circ \al_n)'(0) \vert /2) \subset (\phi \circ \al_n)(B(0, \delta))
\]
for all $n$ large. By (2.2), there is a uniform positive lower bound on $\vert (\phi \circ \al_n)'(0) \vert$, and hence the images $(\phi \circ \al_n)(B(0, \delta))$ contain a ball of fixed size around $\phi \circ \al_n(0)$.

\medskip

For the second step, it suffices to show that for every path $\ga : [0,1] \ra L_p, \ga(0) = p$, there exists a path $\ti \ga: [0,1] \ra \mbb D \sm \ti \al^{-1}(E)$ with $\ti \ga(0) = 0 \in \mbb D$
and $\ti \al \circ \ti \ga = \ga$. We may assume that $\ga$ is $C^1$-smooth. Let $\la_U = g/ \eta_U$. Since $\eta_U$ is assumed to be continuous on $U$, $\la_U$ is continuous on $T \mathcal F_U$, the tangent bundle of $\mathcal F_U$. Also recall that $\la_U$ induces the Poincar\'{e} metric on the leaves and hence
\begin{equation}
k = \sup \left\{ \la_U(\ga(t), \ga'(t)) : t \in [0, 1] \right\}
\end{equation}
is finite since $\ga \subset L_p$ is compact. By \cite{CLN}, there exists a compact neighbourhood $\tilde{V}$ of $\ga([0,1])$ in $L_p$ and a tubular neighbourhood $\pi : \tilde{W} \ra \tilde{V}$ such that if $q \in \pi^{-1}(p)$ is close to $p$, then $\ga$ lifts to a $C^1$-smooth curve $\ga_q : [0,1] \ra L_q$ such that $\pi \circ \ga_q = \ga$ and $\ga_q \ra \ga$ as $q \ra p$. Since $\rho(U, U_n) \ra 0$ and $\tilde{V} \subset \tilde{W}$ are compactly contained in $U$, it follows that $\tilde{W}$ is compactly contained in $U_n$ for $n$ large. Therefore, $\ga \subset L_{p, U_n}$ for $n$ large. 

\medskip

The idea now is to lift $\ga$ using $\pi$ to get paths that are contained in different leaves and then lift them to $\mbb D$ using the uniformizations $\al_n$. To make all this precise, note that as $\al_n(0) \ra \ti \al(0) = p$, the argument using Rouche's theorem given above shows that it is possible to choose a sequence $z_n \in \mbb D$, $z_n \ra 0$ such that $q_n = \al_n(z_n) \in \pi^{-1}(p)$. Note that $q_n \in L_{p_n, U_n}$. Then $\ga$ admits a lift using $\pi$ starting at $q_n$ and contained in $L_{q_n}$. Call this lift $\ga_{q_n}$ and note that $\ga_{q_n} \subset \tilde{W} \subset U_n$ for large $n$ as $\ga \subset \tilde{W}$ and $\ga_{q_n} \ra \ga$. Therefore, there exists a lift $\ti \ga_{q_n} : [0,1] \ra \mbb D$ with $\ti \ga_{q_n}(0) = z_n$ such that $\al_n \circ \ti \ga_{q_n} = \ga_{q_n}$. 

\medskip

\no (d) Let $\mathcal K \subset T\mathcal F_U$ be a compact set. Note that the projection of $\mathcal K$ to $M$ is a compact set that does not intersect the singular set $E$. Since $\rho(U, U_n) \ra 0$, $\mathcal K \subset T \mathcal F_{U_n}$ for all large $n$. We will now compare $\la_U$ and $\la_{U_n}$ on $\mathcal K$.

\medskip

\no {\it Claim:} There exists a constant $C = C(\mathcal K) > 0$ such that
\[
\vert \la_{U_n}(z, v) - \la_U(z, v) \vert \le C
\]
for all $(z, v) \in \mathcal K$ and $n$ large.

\medskip

If this were false, there would exist a sequence of pairs $(w_n, v_n) \in \mathcal K$ such that
\[
\vert \la_{U_n}(w_n, v_n) - \la_U(w_n, v_n) \vert \ge n
\]
after passing to a subsequence. By passing to a further subsequence if needed, $(w_n, v_n) \ra (w_0, v_0) \in \mathcal K$, where $w_0 \notin E$, and by using the continuity of $\la_U$, it can be seen that
\[
\vert \la_{U_n}(w_n, v_n) - \la_U(w_0, v_0) \vert \ge n- 1
\]
for all large $n$. Since $\vert v_n \vert_g, \vert v_0 \vert_g$ are uniformly bounded and $0 < \eta_U(w_0) < \infty$, it follows that $\eta_{U_n}(w_n) \ra 0$. Choose uniformizations $\psi_n : \mbb D \ra L_{w_n, U_n}$ such that $\psi_n(0) = w_n$ and $\eta_{U_n}(w_n) = \vert \psi'_n(0) \vert_g$. Then $\vert \psi'_n(0) \vert_g \ra 0$. But this is a contradiction as the proof of (c) shows. Indeed, since $w_0 \notin E$, there is a uniform positive lower bound on $\vert \psi'_n(0) \vert_g$, and thus the claim holds.

\medskip

\no (e) Coming back to the lifts $\ga_{q_n}$ which are $C^1$-perturbations of $\ga$, it follows from (d) that there is a constant $C > 0$ such that
\[
\sup \left\{ \left \vert \la_{U_n}(\ga_{q_n}(t), \ga'_{q_n}(t)) - \la_U(\ga(t), \ga'(t)) : t \in [0,1] \right \vert \right\} \le C
\]
and hence $\la_{U_n}(\ga_{q_n}(t), \ga'_{q_n}(t)) \le C + k$ for all large $n$ and $t \in [0,1]$. Since unramified holomorphic coverings preserve the infinitesimal hyperbolic metric,
\[
\vert \ti \ga'_{q_n}(t) \vert \le \frac{\vert \ti \ga'_{q_n}(t) \vert}{1 - \vert \ti \ga_{q_n}(t) \vert^2}  \le C + k
\]
and consequently, the maps $t \mapsto \ti \ga_{q_n}(t)$ are equicontinuous. In addition, if $t \in [0,1]$, then $d_h(z_n, \ti \ga_{q_n}(t))$, the hyperbolic distance between $z_n = \ti \ga_{q_n}(0)$ and $\ti \ga_{q_n}(t)$ is at most
\[
\int_0^t \frac{\vert \ti \ga'_{q_n}(t) \vert}{1 - \vert \ti \ga_{q_n}(t) \vert^2} \le  C + k
\]
and since $t$ is arbitrary, this means that the traces of all paths $\ti \ga_{q_n}$ are uniformly compactly contained in $\mbb D$. It follows that, after passing to a subsequence if needed, there is a path $\ti \ga = \lim \ti \ga_{q_n}$ in $\mbb D$ that starts at the origin and satisfies 
\[
\ti \al \circ \ti \ga = \lim \al_n \circ \ti \ga_{q_n} = \lim \ga_{q_n} = \ga.
\]
Since $\ga$ does not intersect $E$, it follows that $\ti \ga$ does not intersect $\ti \al^{-1}(E)$. Thus, $\ti \al : \mbb D \sm \ti \al^{-1}(E) \ra L_p$ is a covering map.

\medskip

\no (f) For $\ep > 0$, note that $\ti \al \circ \psi(\mbb D)$ is compactly contained in $U$, where $\psi : \mbb D \ra \mbb D$ is defined by $\psi(z) = (1- \ep)z$. Therefore, the images $\al_n \circ \psi(\mbb D)$ are also uniformly compactly contained in $U$ for all large $n$. Hence, the maps $\al_n \circ \psi$ are candidates in the family that defines $\eta_U$ and so
\[
(1- \ep) \vert \al_n'(0) \vert_g \le \eta_U(p_n).
\] 
The continuity of $\eta_U$ and the fact that $\al_n \ra \ti \al$ uniformly on compact subsets of $\mbb D$ shows that $(1 - \ep) \vert \ti \al'(0) \vert_g \le \eta_U(p)$. By allowing $\ep \ra 0$, it follows that $\vert \ti \al'(0) \vert_g \le \eta_U(p)$. By choosing a uniformization $\beta : \mbb D \ra L_{p, U}$ with $\be(0) = p$, this inequality can be rephrased as $\vert \ti \al'(0) \vert_g \le \vert \be'(0) \vert_g$.
 
\medskip

\no (g)  Let $\pi : \mbb D \ra \mbb D \sm \ti \al^{-1}(E)$ be a covering with $\pi(0) = 0$. As $\ti \al^{-1}(E) \not= \emptyset$ by assumption, the Schwarz lemma shows that $\vert \pi'(0) \vert < 1$. The following diagram
\begin{center} 
\begin{tikzcd}
\mbb D \arrow{r}{\ti \pi} \arrow[swap]{d}{\pi} & \mbb D \arrow{d}{\beta} \\
\mbb D \sm \ti \al^{-1}(E) \arrow{r}{\ti \al} & L_{p, U}
\end{tikzcd}
\end{center}
commutes, where $\ti \pi : \mbb D \ra \mbb D$ is a lift of $\ti \al \circ \pi : \mbb D \ra L_{p, U}$. Since $\be$ and $\ti \al \circ \pi$ are both uniformizations of $L_{p, U}$ with $\be(0) = \ti \al \circ \pi(0) = p$ , it must be the case that $\vert \ti \pi'(0) \vert = 1$. On the other hand,
\[
\vert \be'(0) \vert_g = \vert (\be \circ \ti \pi)'(0) \vert_g = \vert (\ti \al \circ \pi)'(0) \vert_g = \vert \ti \al'(0) \vert_g \cdot \vert \pi'(0) \vert < \vert \ti \al'(0) \vert_g
\]
and this contradicts the inequality obtained in (d). 

\medskip 

Hence, $\ti \al^{-1}(E) = \emptyset$ which implies that $\ti \al : \mbb D \ra L_{p, U}$ is a uniformization. This completes the proof of Theorem 1.1.


\section{Proof of Theorem 1.2}

\no Recall that $(M, E, \mathcal F)$ is a Riemann surface foliation that is assumed to admit an ultrahyperbolic metric $g$ with curvature bounded above by $-a^2 < 0$. The Ahlfors lemma implies that each leaf $L \subset \mathcal F$ is hyperbolic and $g|L \le a^{-1} \la_L$, where $\la_L$ is the Poincar\'{e} metric on $L$.

\medskip

For a fixed leaf $L \subset \mathcal F$, consider the subfamily $\mathcal O_m(\mbb D, L) \subset \mathcal O_m(\mbb D, \mathcal F)$. For each $f \in \mathcal O_m(\mbb D, L)$, it follows that $f^{\ast}(g |L) \le a^{-1} f^{\ast}(\la_L)$. Hence
\[
f^{\mathcal F, g}_{\mbb D}(0) = \left. \frac{f^{\ast}(g | L)}{\la_{\mbb D}(z) \vert dz \vert}  \right|_{z = 0} \le  a^{-1} \cdot \left. \frac{f^{\ast}(\la_L)}{\la_{\mbb D}(z) \vert dz \vert}  \right|_{z = 0} = a^{-1} \cdot f^L_{\mbb D}(0)
\]
and this means that if $f^{\mathcal F, g}_{\mbb D}(0) \ge \al$, then $f^L_{\mbb D}(0) \ge a \al$. It follows that
\[
\inf \left\{ t(f) : f \in  \mathcal O_m(\mbb D, L), f^{\mathcal F, g}_{\mbb D}(0) \ge \al \right\} \ge \mathcal T_{m, L}(a \al)
\]
where $\mathcal T_{m, L}$ is as in (1.1). The lower bound on $\mathcal T_{m, L}$ in (1.2) holds for all hyperbolic Riemann surfaces and does not require any further properties of $L$. Since each $f \in \mathcal O_m(\mbb D, \mathcal F)$ belongs to $\mathcal O_m(\mbb D, L)$ for some leaf $L \subset \mathcal F$, it follows that
\[
\mathcal T^g_{m, \cal F}(\al) \ge  \log \left( \frac{1 + \tau_m ((a \al)^{1/(m+1)})}{1 - \tau_m ((a \al)^{1/(m+1)})} \right)
\]
for all $\al \in (0, 1/a]$. This completes the proof. 
 

\section{Remarks and Examples}

\no (A)  The results in \cite{N1}, \cite{N2} and \cite{NM} describe a paradigm that can be summarized thus. Let $g$ be a hermitian metric on $M\sm E$ whose restriction to each leaf $L \subset \mathcal F$ has Gaussian curvature bounded from above by a uniform negative constant. Then each leaf is hyperbolic. Furthermore, if $g$ is complete, then under suitable hypotheses on the singular locus $E$, the family of uniformizations of the leaves of $\mathcal F$ is normal, i.e., a family of uniformizations has a subsequence that either converges to a uniformization of the limiting leaf or to a constant map with values in $E$. A direct consequence of this is the continuity of $\eta$. As remarked before the statement of Theorem 1.1, it is possible to consider localized versions of $\eta$ with respect to a given neighbourhood in $X$, and Theorem 1.1 shows that even these localized versions of $\eta$ possess continuity properties as the neighbourhoods vary in the Hausdorff sense.   

\medskip

Several examples on the growth of $\eta$ near the singular set $E$ can be found in \cite{N1}, \cite{N2} and \cite{NM}. They are based on explicit local constructions of metrics that reflect the intrinsic geometry of the vector fields defining the foliation $\mathcal F$. It turns out that these examples are also useful in illustrating the various possibilities that can arise in Theorem 1.2 above. We collate them below keeping this connection in mind. In what follows, $X$ is a germ of a holomorphic vector field near the origin in $\mbb C^n$. We will work with a representative of $X$ that is defined in the ball $\mbb B^n(0, \rho)$  of radius $\rho < 1$ and the standing assumption is that the origin is an isolated singularity for $X$. Any additional hypotheses will be mentioned explicitly in the various examples given below. The resulting foliation will be denoted by $\mathcal F$.

\medskip

\no {\it Example 1:} Consider the restriction of the smooth metric
\[
g(z) = \frac{2 \; \vert dz \vert}{1 - \vert z \vert^2} 
\]
to the leaves of $\mathcal F$. Let $Z = Z(T)$ be the solution of
\[
\dot{Z} = dZ/dT = X(Z(T)), \; Z(0) = z \in \mbb B^n(0, \rho) \sm \{0\}
\]
where $Z$ is defined in a neighbourhood of the origin in the plane. Since $\dot{Z} = X(Z)$, a computation shows that the curvature of the pull back metric $Z^{\ast}g$ at $T = 0$ is
\[
K_g(z) = (\vert z \vert^2 - 1) - \vert X(z) \vert^{-2} \left | \left \langle z, X(z) \right \rangle \right |^2 - \frac{1}{2}(1 - \vert z \vert^2)^2 \vert X(z) \vert^{-6} \left \vert X(z) \wedge DX(z)(X(z)) \right \vert^2
\]
where the angle-brackets denote the standard hermitian inner product on $\mbb C^n$ and $\vert z \wedge w \vert^2 = \vert z \vert^2 \vert w \vert^2 - \vert \langle z, w \rangle \vert^2$. Since $\vert z \vert < \rho < 1$,
\[
K_g(z) \le \vert z \vert^2  - 1 < \rho^2 - 1 < 0
\]
and hence $g$ is ultrahyperbolic of curvature bounded by $-a^2 < 0$ where $a^2 = 1 - \rho^2 > 0$. Theorem 1.2 then shows that 
\[
\mathcal T^g_{m, \mathcal F_{\rho}}(\al) \ge  \log \left( \frac{1 + \tau_m ((a \al)^{1/(m+1)})}{1 - \tau_m ((a \al)^{1/(m+1)})} \right)
\]
for all $\al \in (0, 1/a]$.

\medskip

\no {\it Example 2:} Suppose that the origin is a non-degenerate singularity of $X$ in the sense that there exists $C > 1$ such that
$C^{-1}\vert z \vert \le \vert X(z) \vert \le C \vert z \vert$ for $z$ near the origin. Let $h$ be a smooth hermitian metric on $\mbb B^n(0, \rho)$ and let $u$ be a continuous plurisubharmonic function on $\mbb B^n(0, \rho)$ which is $C^{\infty}$-smooth away from the origin. Consider the restriction of the metric
\[
g = \frac{e^u}{(\log \vert z \vert)^2} \cdot \frac{h}{h(X)}
\]
which is smooth away from the origin, to the leaves of $\mathcal F$. If $Z = Z(T)$ is a solution of
\[
\dot{Z} = dZ/dT = X(Z(T)), \; Z(0) = z \in \mbb B^n(0, \rho) \sm \{0\}
\]
near $T = 0$, it can be seen that the pull back metric $Z^{\ast}g = f(T) \vert dT \vert^2$ where 
\[
f(T) = \frac{e^{u(Z(T))}}{(\log \vert Z(T) \vert)^2}.
\]
and its curvature $K_g$ satisfies 
\[
-\frac{1}{2}K_g(z) = \frac{1}{f(0)} \pa \overline \pa (u(Z(T)) + \log(\log^{-2}(\vert Z(T) \vert))) \big|_{T=0}.
\]
The plurisubharmonicity of $u$ implies that the first term above is non-negative and hence
\[
-\frac{1}{2}K_g(z) \ge e^{-u(z)}(\log \vert z \vert)^2 \pa \overline \pa \log(\log^{-2}(\vert Z(T) \vert)) \big|_{T=0}.
\]
Now,
\[
\pa \overline \pa \log(\log^{-2}(\vert Z(T) \vert)) \big|_{T=0} = \frac{1}{\vert z \vert^4 \vert \log (\vert z \vert) \vert} \left( \vert z \vert^2 \vert X(z) \vert^2 - \vert \langle X(z), z \rangle \vert^2  \right) + \frac{\vert \langle X(z), z \rangle \vert^2}{2 \vert z \vert^4 \log^2(\vert z \vert)}
\]
as a computation shows and the Cauchy-Schwarz inequality then implies that $K_g(z) < 0$ away from the origin. Note that 
\begin{equation}
-\frac{1}{2}K_g(z) \ge \frac{\vert \log (\vert z \vert) \vert}{e^{u(z)}} \left( \frac{\vert X(z) \vert^2}{\vert z \vert^2} - \frac{\vert \langle X(z), z \rangle \vert^2}{\vert z \vert^4}  \right) + \frac{\vert \langle X(z), z \rangle \vert^2}{2 \vert z \vert^4 e^{u(z)}}.
\end{equation}
It follows that if $\vert z \vert <  \rho < 1/\sqrt{e}$, then $\vert \log(\vert z \vert) \vert > 1/2$ and hence
\[
-K_g(z) \ge \frac{\vert X(z) \vert^2}{\vert z \vert^2 e^{u(z)}}
\] 
as can be seen by combining the second and third terms in (4.1). Finally, using the assumption that the origin is a non-degenerate singularity, it follows that 
\[
\limsup_{z \ra 0} K_g(z) \le  -C^{-2}/e^{u(0)}
\]
and hence that $K_g(z)$ is bounded above by a negative constant, say $-a^2$, in a neighbourhood of the origin in case $u(0) > -\infty$. The same conclusion holds even if $u(0) = -\infty$. These calculations can be found in \cite{N1} under the hypotheses that $u$ is pluriharmonic near the origin. As indicated above, they hold even when $u$ is a continuous plurisubharmonic function.

\medskip

\no {\it Example 3:} Suppose that $\mathcal F$ admits $n$ smooth invariant hypersurfaces that meet in general position at the origin. In this case, there is a coordinate system in which the  invariant hypersurfaces are given by $\{z_j = 0\}$, $1 \le j \le n$ and 
\[
X(z) = z_1^{k_1} f_1(z) \; \pa / \pa z_1 + z_2^{k_2} f_2(z) \; \pa / \pa z_2 + \ldots + z_n^{k_n} f_n(z) \; \pa / \pa z_n
\]
where the $f_i$'s are holomorphic near the origin and $k_i \ge 1$ are integers. For $z_j \not= 0$, set 
\[
\psi_j(z) = \frac{\vert z_j \vert^{2k_j-2} \vert f_j(z) \vert^2}{\log^2(\vert z_j \vert)}
\]
and note that $\psi_j$ extends to a continuous function on $\mbb B^n(0, \rho) \sm \{0\}$ by setting $\psi_j(z) = 0$ if $z_j = 0$. Then $\psi(z) = \max \{ \psi_j(z) : 1 \le j \le n\}$ is a positive continuous function away from the origin. As shown in \cite{NM}, the metric
\[
g(z) = \frac{\psi(z)}{\vert X(z) \vert^2} \; \vert dz \vert^2
\]
is ultrahyperbolic of curvature bounded above by $-1$. 
 
\medskip

All these examples show that it is possible to define smooth or ultrahyperbolic metrics away from the origin with curvatures bounded above by a negative constant. The uniform negative upper bound for the curvature is determined by the metric which in turn reflects the properties of the underlying holomorphic vector field. Theorem 1.2 shows that in each of these examples, 
the foliated domain Bloch constant
\[
\mathcal T^g_{m, \mathcal F_{\rho}}(\al) \ge  \log \left( \frac{1 + \tau_m ((a \al)^{1/(m+1)})}{1 - \tau_m ((a \al)^{1/(m+1)})} \right)
\]
for all $\al \in (0, 1/a]$. Here, the superscript $g$ refers to the metrics in the examples  above, and $a > 0$ controls the curvature of the leaves near the origin in the corresponding cases.

\medskip

\no (B) To propose a minor variation of Minda's construction of the domain Bloch constant $\mathcal T_{m. S}(\al)$ in (1.1), let $\al : S \ra [0,1]$ be a continuous function on the hyperbolic Riemann surface $S$. For $x \in S$, define
\[
\mathcal T^{\al}_{m, S}(x) = \inf \left\{ t(f) : f \in \mathcal O_m(\mbb D, S), f(0) = x, f^S_{\mbb D}(0) \ge \al(x) \right\}
\] 
which allows for a variable lower bound for the distortion function. Since the family that defines $\mathcal T^{\al}_{m, S}(x)$ is a subset of that which defines (1.1), it follows that
\[
\mathcal T^{\al}_{m, S}(x) \ge \mathcal T_{m, S}(\al(x)) \ge \log \left( \frac{1 + \tau_m ((\al(x))^{1/(m+1)})}{1 - \tau_m ((\al(x))^{1/(m+1)})} \right).
\]
With this in place, it is straightforward to formulate analogs of Theorem 1.2 in which we consider a continuous function $\al : \mathcal F \ra [0,1]$ and use it to allow variable lower bounds for the distortions of maps $f \in \mathcal O_m(\mbb D, \mathcal F)$ and finally obtain a uniform lower bound for the ball on which every such $f$ is injective.  


\end{document}